\documentclass[11pt,oneside]{amsart}
\usepackage{amssymb}

\usepackage{cite}

\DeclareFontFamily{OT1}{rsfs}{}
\DeclareFontShape{OT1}{rsfs}{m}{n}{ <-7> rsfs5 <7-10> rsfs7 <10-> rsfs10}{} 
\DeclareMathAlphabet{\mathscr}{OT1}{rsfs}{m}{n}

\newcommand{\eq}[1]{\eqref{#1}}

\newcommand{\bel}[1]{\begin{equation}\label{#1}}
\newcommand{\beal}[1]{\begin{eqnarray}\label{#1}}
\newcommand{\beadl}[1]{\begin{deqarr}\label{#1}}
\newcommand{\eeadl}[1]{\arrlabel{#1}\end{deqarr}}
\newcommand{\eeal}[1]{\label{#1}\end{eqnarray}}
\newcommand{\eead}[1]{\end{deqarr}}
\newcommand{\eea}{\end{eqnarray}}
\newcommand{\eeaa}{\end{eqnarray*}}

\newcommand{\be}{\begin{equation}}
\newcommand{\ee}{\end{equation}}

\DeclareFontFamily{OT1}{rsfs}{}
\DeclareFontShape{OT1}{rsfs}{m}{n}{ <-7> rsfs5 <7-10> rsfs7 <10->
rsfs10}{} \DeclareMathAlphabet{\mycal}{OT1}{rsfs}{m}{n}

\newcounter{mnotecount}[section]

\newcommand{\rmnote}[1]{}

\newcommand{\Ric}{\operatorname{Ric}}

\def\mysavedown#1{\edef\mysubs{\mysubs#1}}
\def\mysaveup#1{\edef\mysups{\mysups#1}}
\def\mydown#1{{\mytensor}_{\vphantom{\mysubs}#1}}
\def\myup#1{{\mytensor}^{\vphantom{\mysups}#1}}
\def\tensor#1#2{
  #1
  \def\mytensor{\vphantom{#1}}
  \def\mysubs{\relax}
  \def\mysups{\relax}
  \let\down=\mysavedown
  \let\up=\mysaveup
  #2
  \let\down=\mydown
  \let\up=\myup
  #2
  }

\newcommand{\letters}
  {\renewcommand{\theenumi}{\alph{enumi}}
   \renewcommand{\labelenumi}{(\theenumi)}}
\newcommand{\romanletters}
  {\renewcommand{\theenumi}{\roman{enumi}}
   \renewcommand{\labelenumi}{(\theenumi)}}

\newcommand{\pd}[1]{\frac{\del}{\del #1}}

\newcommand{\Tr}{\operatorname{Tr}}
\newcommand{\Id}{\operatorname{Id}}

\newcommand{\R}{\mathbb R}

\renewcommand{\to}{\rightarrow}

\newcommand{\Ccirc}{\mathring{C}{}}

\renewcommand{\centerdot}{\mathbin{\text{\protect\raisebox{-.3ex}[1ex][0ex]{\Large{$\cdot$}}}}}
\newcommand{\cross}{\mathbin{\times}}

\newcommand{\Del}{\nabla}

\renewcommand{\exp}{\operatorname{exp}}

\DeclareMathOperator{\vectorspan}{span}
\newcommand{\half}{{\tfrac12}}
\newcommand{\mapsinto}{\mathrel{\hookrightarrow}}
\renewcommand{\phi}{\varphi}
\renewcommand{\epsilon}{\varepsilon}
\renewcommand{\hat}{\widehat}
\let\scr=\mathscr

\def\crn#1#2{{\vcenter{\vbox{
        \hbox{\kern#2pt \vrule width.#2pt height#1pt
           }
          \hrule height.#2pt}}}}

\newcommand{\Kbar}{{\smash[t]{\overline K}}}
\newcommand{\tw}{\widetilde}

\newcommand{\<}{\langle}
\renewcommand{\>}{\rangle}
\newcommand{\del}{\partial}
\renewcommand{\H}{\mathbb H}

\newcommand{\grad}{\operatorname{grad}}

\renewcommand{\hbar}{{\overline h}}

\renewcommand{\)}{\textup{)}}

\newcommand{\pre}[2]{{{\vphantom{#2}}^{#1}}\kern-.2ex{#2}}

\theoremstyle{plain}
\newtheorem{theorem}{Theorem}[section]
\newtheorem{lemma}[theorem]{Lemma}

\newtheorem{bigtheorem}{Theorem}[section]

\theoremstyle{definition}

\theoremstyle{remark}

\numberwithin{equation}{section}

\date{\today}

\begin{document}

\title[Boundary regularity]
{Boundary regularity of conformally compact
Einstein metrics}

\author[Chru\'sciel]{Piotr T. Chru\'sciel}
\address{Piotr T. Chru\'sciel,
D\'epartement de math\'ematiques,
Facult\'e des Sciences,
Parc de Grandmont,
F37200 Tours, France}
\email{Piotr.Chrusciel@lmpt.univ-tours.fr}
\urladdr{http://www.phys.univ-tours.fr/\~{}piotr}
\thanks{First author partially supported by Polish Research Council grant KBN 2 P03B 073 24;
second author partially supported by the ACI program of the French Ministry of
Research.}
\author[Delay]{Erwann Delay}
\address{Erwann Delay,
D\'epartement de math\'ematiques,
Facult\'e des Sciences, rue Louis Pasteur, F84000 Avignon,
France}
\email{Erwan.Delay@univ-avignon.fr}
\urladdr{http://www.phys.univ-tours.fr/\~{}delay}
\author[Lee]{John M. Lee}
\address{John M. Lee,
University of Washington,
Department of Mathematics,
Box 354350,
Seattle, WA 98195-4350}
\email{lee@math.washington.edu}
\urladdr{http://www.math.washington.edu/\~{}lee}
\author[Skinner]{Dale N. Skinner}
\address{Dale N. Skinner,
Milliman USA,
1301 Fifth Avenue, Suite 3800,
Seattle, WA 98101}
\email{dales82@mindspring.com}
\urladdr{http://www.daleskinner.com}

\begin{abstract}
We show that $C^{2}$ conformally compact Riemannian Einstein
metrics have conformal compactifications that are smooth up to the
boundary in dimension $3$ and all
even dimensions, and polyhomogeneous in odd
dimensions greater than $3$.
\end{abstract}

\maketitle

\section{Introduction}\label{section:intro}
Suppose $\overline M$ is a smooth, compact
manifold with boundary, and let $M$ denote its interior and $\del
M$ its boundary. (By ``smooth,'' we always mean $C^\infty$.) A
Riemannian metric $g$ on $M$ is said to be {\it conformally
compact} if for some smooth defining function $\rho$ for $\del M$
in $\overline M$, $\rho^2 g$ extends by continuity to a Riemannian
metric (of class at least $C^0$) on $\overline M$. The rescaled
metric $\overline g = \rho^2 g$ is called a {\it conformal
compactification} of $g$. If for some (hence any) smooth defining
function $\rho$, $\overline g$ is in $C^{k}(\overline M)$ or
$C^{k,\lambda}(\overline M)$, then we say $g$ is conformally
compact of class $C^k$ or $C^{k,\lambda}$, respectively.

If $g$ is conformally compact, the restriction of $\overline g = \rho^2 g$
to $\del M$ is a Riemannian metric on $\del M$, whose conformal class
is determined by $g$, independently of the choice of defining function $\rho$.
This conformal class is called the
{\it conformal infinity} of $g$.

Several important existence and uniqueness results
\cite{mand2,mand1,biquard,GL,Lee:fredholm}
concerning conformally compact Riemannian Einstein metrics
have been established recently.  For many applications to physics and
geometry, it turns out to be of great importance to understand the 
asymptotic behaviour of the resulting metrics near the boundary.  This
question has been addressed by Michael Anderson
\cite{mand1}, who proved that if $g$ is a $4$-dimensional
conformally compact Einstein metric with smooth conformal
infinity, then the conformal compactification of $g$ is smooth up
to the boundary in suitable coordinates.  
It has long been conjectured that in higher dimensions, 
conformally compact Einstein metrics with
smooth conformal infinities should have infinite-order asymptotic expansions 
in terms of $\rho$ and $\log\rho$.  The purpose of this paper is to 
confirm that conjecture.

The choice of special coordinates in Anderson's result cannot be
dispensed with. Because the Einstein equation is invariant under
diffeomorphisms, we cannot expect that the conformal
compactification of an arbitrary conformally compact Einstein
metric will necessarily 
have optimal regularity for all $C^\infty$ structures on $\overline M$.
For example, suppose $g$ is an Einstein metric
on $M$ with a smooth conformal compactification, and let $\Psi \colon
\overline M\to \overline M$ be a
homeomorphism that restricts to the
identity map of $\del M$ and to a diffeomorphism from $M$
to itself.
Then $\Psi ^*g$ will still
be Einstein with the same conformal infinity, but its conformal
compactification $\rho^2 \Psi ^*g$ may no longer be smooth. Thus the best
one might hope for is that an arbitrary conformally compact
Einstein metric can be made smoothly conformally compact after
pulling back by an appropriate diffeomorphism. Even this is not true
in general, because Fefferman and Graham showed in \cite{FG} that
there is an obstruction to smoothness in odd dimensions.

Since Einstein metrics are always smooth (in fact, real-analytic) in suitable coordinates
in the interior, only regularity at the boundary is at issue.
For that reason, instead of assuming
that $\overline M$ is compact, we will 
assume only that it
has a compact boundary component $Y$, and restrict our attention
to a collar neighborhood
of $Y$ in $\overline M$, which we may assume without loss of generality 
is diffeomorphic to $Y\cross[0,1)$.
Throughout this paper, then, $Y$ will be an arbitrary smooth, connected, compact, 
$n$-dimensional manifold without boundary, and we make the following identifications:
\begin{displaymath}
\overline M=Y\cross[0,1), \quad 
M = Y\cross(0,1), \quad
\del M = Y\cross\{0\}.
\end{displaymath}
Let $\rho\colon \overline M\to [0,1)$ denote the projection onto the $[0,1)$ 
factor; it is a smooth
defining function for $\del M$ in $\overline M$.
For $0<R<1$, we define%
\begin{displaymath}
M_R =  Y\cross(0,R],\quad
\overline M_R =  Y\cross[0,R].
\end{displaymath}
In this context, we extend the definition of conformally compact metrics by 
saying that a Riemannian metric $g$ on $M$ or $M_R$ is conformally compact
if $\rho^2 g$ extends to a continuous metric on $\overline M$ or $\overline M_R$, 
respectively.
A continuous map
$\Psi \colon
\overline M_R\to \overline M$ (for some $R$) that restricts to the 
identity map of $\del M$ and to a diffeomorphism from $M_R$
to its image will be called a
{\it collar diffeomorphism}.
If $\Psi $ and its inverse are of class $C^k$ (or $C^{k,\lambda}$) up to the boundary, we will
call it a $C^k$ (resp., $C^{k,\lambda}$) collar diffeomorphism.

The object of this paper is to prove that
the following regularity holds:

\begin{bigtheorem}\label{thm:main}
Let
$g$ be a Riemannian metric on $M$.
Suppose that $\dim M = n+1\ge 3$;
$g$ is Einstein with $\Ric(g) = - n g$; $g$ is
conformally compact of class $C^2$; and the
representative $\gamma = \rho^2 g|_{\del M}$ of the conformal
infinity of $g$ is smooth.
Let $\tw\gamma$ be any smooth representative of the conformal class
$[\gamma]$.
Then for any $0<\lambda<1$, there exists $R>0$ and
a $C^{1,\lambda}$
collar diffeomorphism $\Phi \colon \overline M_R\to \overline M$
such that
$\Phi ^*g$
can be written in the form
\begin{equation}\label{eq:metric-in-Fermi-coords}
\Phi ^*g = \rho^{-2} (d\rho^2 + G(\rho)),
\end{equation}
where $\{G(\rho):0<\rho\le R\}$ is a one-parameter family of smooth
Riemannian metrics on $Y$,  $d\rho^2 + G(\rho)$
has a continuous
extension to $\overline M_R$ with
$G(0) = \tw\gamma$, and has the following regularity:
\begin{enumerate}\letters
\item
If $\dim M$ is even or equal to $3$, then $d\rho^2 + G(\rho)$ extends 
smoothly to
$\overline M_R$, so $\Phi^*g$ is conformally compact of class $C^\infty$.
\item
If $\dim M$ is odd and greater than $3$, then
$G$ can be written in the form
\begin{displaymath}
G(\rho) = \phi(\rho,\rho^n\log \rho),
\end{displaymath}
with $\phi(\rho,z)$ a two-parameter family of Riemannian metrics on $Y$
that is smooth in all of its arguments as a function
on $Y\cross [0,R]\cross [R^n \log R,0]$.
Furthermore, $\Phi^*g$ is smoothly conformally compact if and only
if $\del_z\phi(0,0)$ vanishes identically on $\del M$.
\end{enumerate}
\end{bigtheorem}

\noindent{\it Remark.}  The symmetric $2$-tensor field
$\del_z\phi(0,0)$ along $\del M$ 
can be determined in principle from local
computations involving the conformal class $[\gamma]$
(cf.\/~\cite{FG,FG2003,HG}).  In fact, this tensor field is a constant
multiple of the {\it ambient obstruction tensor} defined
in \cite{HG}, whose vanishing
is a necessary condition for the existence of a 
smoothly conformally compact Einstein metric on $M$ with
$[\gamma]$
as conformal infinity.
Explicit formulae in low dimensions can be found in~\cite{deHaro:2000xn}.
In odd dimensions, it is shown in~\cite{KichenassamyFG} that for
any analytic $[\gamma]$ and $\partial_z \varphi(0,0)$ there exists
a unique Einstein metric as above, defined on some neighborhood of
the conformal boundary.

The main idea of the proof is to use the harmonic map equation to put
$g$ into a gauge in which it satisfies an elliptic equation, and
then apply the polyhomogeneity results of \cite{AndChDiss}. The
proof consists of four steps.  First, we construct a preliminary
collar diffeomorphism that makes $\rho^2 g$ coincide to second order
along $\del M$ with a smooth product metric $\overline h$.
Second, applying the inverse
function theorem to the harmonic map equation, we show that there
exists a collar diffeomorphism $H\colon \overline M_R\to \overline M_R$
that is harmonic in $M_R$,
thought of as a map from $(M_R,g)$ to
$(M_R,h)$, where $h = \rho^{-2}\overline h$.  It follows that the
metric $\tw g = (H^{-1})^* g$ satisfies the following
``gauge-broken Einstein equation'' near $\del M$:%
\begin{equation}\label{eq:gauge-broken-Einstein}
Q(\tw g,h) := \Ric({\tw g}) + n \tw g - \delta_{\tw g}^*(\Delta_{\tw gh}(\Id)) = 0
\end{equation}
(see, e.g.,\/~\cite{Lee:fredholm}), where $\Delta_{\tw gh}$ is the
harmonic map Laplacian. The third step is to show that solutions
to \eqref{eq:gauge-broken-Einstein} satisfy the hypotheses of
\cite[Theorem 5.1.1]{AndChDiss} and therefore are polyhomogeneous
(i.e., have asymptotic expansions in powers of $\rho$ and
$\log\rho$).  The last step is to use a special defining function
and Fermi coordinates near the boundary to put the metric into the
form \eqref{eq:metric-in-Fermi-coords}.

In this paper, we have addressed the regularity issue only for 
{\it smooth} conformal infinities, primarily because the 
polyhomogeneity results of \cite{AndChDiss} are proved only in
that context.  
By chasing the losses of differentiability that occur at
various steps of our construction and that of 
\cite{AndChDiss}, 
one can almost certainly prove
that for any $N$
there exists $n(N)$ such that the solution will have a partial
polyhomogeneous expansion with $C^N$ coefficients.  Working out
$n(N)$ would be a straightforward but extremely tedious exercise, which
we have not attempted to do.  

The first and third authors would like to acknowledge useful discussions
with Michael Anderson and Robin Graham.

\section{Weighted H\"older spaces}\label{sec:def}

Throughout most of this paper, we will use the notations and conventions
of \cite{Lee:fredholm}.
We define $M_R$ and $\overline M_R$ as in the introduction.
We assume throughout that $\dim M = n+1\ge 3$.
Any smooth local coordinates $\theta = (\theta^1,\dots,\theta^n)$
on an open set $U\subset Y$ yield smooth coordinates
$(\theta,\rho) = (\theta^1,\dots,\theta^n,\rho)$
on the open subset  $\Omega = U\cross[0,1)\subset \overline M$.
Choose finitely many such charts $(U_i)$ to cover $Y$,
with each set $U_i$ chosen so that
the coordinate functions
extend smoothly to a neighborhood of $\overline U_i$ in $ Y$.
The resulting coordinates on
$\Omega_i = U_i\cross [0,1)\subset \overline M$
will be called {\it background coordinates} for $\overline M$.

Let $B_1,B_2$ be fixed precompact open coordinate balls in 
the upper half-space $\H^{n+1}=\{(x,y)= (x^1,\dots,x^n,y): y>0\}$,
with $(0,\dots,0,1) \in B_1 \subset \overline B_1\subset B_2$.
We will use the summation convention, with Greek indices generally understood
to run from $1$ to $n$, and
Roman indices to run from $1$ to $n+1$; sometimes it will be
convenient to denote $\rho$ by $\theta^{n+1}$ and $y$ by $x^{n+1}$.
Suppose $p\in M_R$,
and let $(\theta_0,\rho_0)$ be the coordinate representation
of $p$ in some fixed background chart.
If $p$ is sufficiently close to $\del M$, we can
define a diffeomorphism $\Phi_p\colon \overline B_2 \to M$  by
\begin{displaymath}
(\theta,\rho) = \Phi_p(x,y) = (\theta_0 + \rho_0 x, \rho_0 y).
\end{displaymath}
As is shown in \cite{Lee:fredholm}, for $R$ sufficiently small,
there exists a countable set of points $\{p_i\}\subset M_R$ such
that the sets $\left\{\Phi_{p_i}(\overline B_2)\right\}$ form a
uniformly locally finite covering of $M_R$, and the sets
$\{\Phi_{p_i}(B_1)\}$ still cover $M_R$. For each such map, set
$\Phi_i = \Phi_{p_i}$, $V_1(p_i) = \Phi_{i}(B_1)$ and $V_2(p_i) =
\Phi_{i}(B_2)$. Then for each $i$,
$\big(V_2(p_i),\Phi_{i}^{-1}\big)$ is a coordinate chart on $M_R$,
called a {\it M\"obius chart}; the corresponding coordinates
$(x,y)$ will be called {\it M\"obius coordinates}. 
It is shown in
\cite[Lemma 2.1]{Lee:fredholm} that if $g$ is any $C^{k,\lambda}$
conformally compact metric on $M$, the pulled-back metrics
$\Phi_i^*g$ are all uniformly $C^{k,\lambda}$ equivalent
to the hyperbolic metric $y^{-2}\sum_j (dx^j)^2$ on $\overline B_2$.
By
compactness, they are also uniformly $C^{k,\lambda}$ equivalent on
$\overline B_2$ to the Euclidean metric $\sum_j (dx^j)^2$.

We will be working in weighted H\"older spaces whose norms
reflect the intrinsic geometry of a conformally compact metric.
Before introducing them,
let us record some elementary facts about H\"older spaces on
subsets of $\R^m$. If
$U\subset \R^m$ is a precompact open subset, $k,p$ are
nonnegative integers, and $\lambda\in[0,1)$,
we denote by
$C^{k,\lambda}(\overline U;\R^p)$ the standard H\"older
space of functions from $\overline U$ to $\R^p$,
and we denote the usual H\"older norm on this space by
$\|\centerdot\|_{k,\lambda;\overline U}$.
If $V\subset \R^p$ is any open set,
$C^{k,\lambda}(\overline U;V)$ denotes the open subset of
$C^{k,\lambda}(\overline U;\R^p)$ consisting of maps that take
their values in $V$.   

\begin{lemma}\label{lemma:holder-facts}
Let $U\subset \R^m$ and $V\subset \R^p$ be convex, precompact open sets,
let $k$ be a nonnegative integer, and let $\lambda\in[0,1)$.
Given $f\in C^{k+2}(\overline V;\R)$ and
$u_0\in C^{k,\lambda}(\overline U;V)$ with $k\ge 1$, there exists
$\delta>0$ and a constant $C = C(U,V,k,\lambda,f,u_0,\delta)$ such that
the following estimate holds for all $u\in C^{k,\lambda}(\overline U;V)$
with $\|u-u_0\|_{k,\lambda;\overline U}\le\delta$:
\begin{displaymath}
\|f\circ u - f\circ u_0\|_{k,\lambda;\overline U} \le C
\|u-u_0\|_{k,\lambda;\overline U}.
\end{displaymath}
\end{lemma}

\begin{proof}
This follows easily from the fact that 
composition $u\mapsto f\circ u$ defines a $C^1$ map from $C^{k,\lambda}(\overline U;V)$
to $C^{k,\lambda}(\overline U;\R)$ (cf.~ \cite{llave-obaya}).
\end{proof}

Now we proceed to define our weighted H\"older spaces
on $M_R$.
Let $E\to M$ be any tensor bundle, and let $R>0$ be chosen so that
$M_R$ is covered by M\"obius charts as above.
For an integer $k\ge 0$ and $\lambda\in(0,1)$,
we define the intrinsic H\"older space
$C^{k,\lambda}(M_R;E)$ to be the set of
locally $C^{k,\lambda}$ sections of $E$ over $M_R$ whose
component functions in M\"obius coordinates satisfy a uniform $C^{k,\lambda}$ bound,
with norm
\begin{displaymath}
\|u\|_{k,\lambda} := \sup_i \|\Phi_{i}^* u\|_{k,\lambda;\overline B_2},
\end{displaymath}
where the supremum is over the
countable collection of M\"obius charts described above.
Weighted versions of these H\"older spaces are defined by setting
$C^{k,\lambda}_\delta (M_R;E) = \rho^\delta  C^{k,\lambda}(M_R;E)$, with the norm
$\|u\|_{k,\lambda,\delta} := \|\rho^{-\delta}u\|_{k,\lambda}$.
It is shown in \cite[Lemma 3.5]{Lee:fredholm} that
this norm is equivalent to
\begin{displaymath}
\|u\|_{k,\lambda,\delta} \sim \sup_i \rho(p_i)^{-\delta}\|\Phi_{i}^* u\|_{k,\lambda;\overline B_j}
\end{displaymath}
for either $j=1$ or $j=2$.

\section{A preliminary normalization}

Suppose that $g$ satisfies the hypotheses of Theorem \ref{thm:main}.
Let 
$\tw\gamma$ be an arbitrary smooth
representative of the conformal class $[\gamma]$ on $Y$.
Define a smooth product metric $\overline h$ on
$\overline M =  Y\cross [0,1)$ 
by
\begin{displaymath}
\overline h = d\rho^2 + \tw\gamma.
\end{displaymath}
Let
$h = \rho^{-2} \overline h$, which is smoothly conformally compact
and has $[\tw\gamma]$ as conformal infinity. 

The goal of this section 
is to show that we can modify $g$ by a 
collar diffeomorphism so that it agrees with $h$ to second order along $\del M$.
This lemma
requires somewhat more work than might be expected
because we are only assuming that $g$ has a $C^2$ conformal compactification,
and we need to make this normalization without losing any smoothness.

\begin{lemma}\label{lemma:bdry-adjusting-diffeo}
Let $g$ and $h$ be as above.
For any 
sufficiently small $R>0$, there exists a $C^3$ collar
diffeomorphism $G\colon \overline M_R\to \overline M$ that
satisfies
$\rho^2 G^* g = \rho^2 h + O(\rho^2)$ in
any background coordinates.
\end{lemma}

\begin{proof}
After replacing $g$ by $G_0^*g$, 
where $G_0(x,\rho) = (x,\rho/f(x))$, we may as well 
assume that $\rho^2 g|_{T\del M} = \tw\gamma$.
We will begin by showing that there exists
a $C^3$ defining function $r$
satisfying
$|dr|^2_{r^2 g} = 1+O(\rho^2)$.
Let 
$\overline g = \rho^2 g$, which is a
$C^2$ Riemannian metric on $\overline M_R$.  
Because $\Ric(g) = -ng$, it follows that
$|d\rho|_{\overline g}^2 =1$ along $\del M$
(cf.\/~\cite[p.\ 192]{GL}).
By Taylor's
theorem, therefore, there is a function $b\in C^1(\overline M_R)$ such that
\begin{displaymath}
|d\rho|_{\overline g}^2 =1 + b\rho.
\end{displaymath}
By Corollary 3.3.2 of \cite{AndChDiss},
there exists a real-valued
function $r\in C^3(\overline M_R)\cap C^\infty(M_R)$
such that
\begin{align*}
r|_{\del M} &= 0,\\
\del_\rho r|_{\del M} &=  1,\\
\del_\rho^2 r|_{\del M} &=  - b|_{\del M}.
\end{align*}
Then $r$ is a $C^3$ defining function for $\del M$, which satisfies
$r = \rho - \half b\rho^2+O(\rho^3)$, $dr = (1-b\rho)d\rho + O(\rho^2)$.
Therefore,
\begin{equation}\label{eq:dr2}
\begin{aligned}
|dr|^2_{r^2 g} &=
|dr|^2_{(r/\rho)^2\overline g}\\
&= (r/\rho)^{-2} |dr|^2_{\overline g}\\
&= (1- \half b\rho)^{-2} (1-b\rho)^2|d\rho|^2_{\overline g} + O(\rho^2)\\
&= 1+O(\rho^2).
\end{aligned}
\end{equation}

Let $P$ denote the gradient of $r$ with respect to the metric $r^2 g$;
thus $P$ is a $C^2$
vector field on $\overline M_R$.  Because $P\rho = \<dr,d\rho\>_{r^2 g}
= \<dr,dr+ O(\rho)\>_{r^2 g} = 1+O(\rho)$,
we can write $P = \del/\del\rho+Q$, where $Q$
is a $C^2$ vector field on $\overline M_R$ that is tangent to $\del M$.
Choose a
smooth embedding $X = (X^1,\dots,X^N)\colon  Y\mapsinto \R^N$ 
into some Euclidean space,
and denote by the same symbol 
the extension of each coordinate function $X^A$ to 
$\overline M =  Y\cross [0,1)$, chosen to be 
constant along the $[0,1)$ factor. 
By \cite[Cor.~3.3.2]{AndChDiss} again, for each
$A=1,\dots,N$,
there is a $C^3$ function 
$\tw X^A\colon \overline M\to \R$, smooth
in $M$, satisfying
\begin{align*}
\tw X^A|_{\del M}&= X^A|_{\del M},\\
\del_\rho \tw X^A|_{\del M} &= - QX^A|_{\del M},\\
\del_\rho^2 \tw X^A|_{\del M} &= - \del_\rho QX^A|_{\del M}.
\end{align*}
It follows that $P\tw X^A = O(\rho^2)$.
The $C^3$ map
$\tw X\colon \overline M\to \R^N$
whose coordinate functions are $(\tw X^1,\dots,\tw X^N)$
thus satisfies $\tw X_*P = O(\rho^2)$.

We wish to use $\tw X$ and $r$ to construct a collar diffeomorphism of $\overline M_R$.
However, $\tw X$ might not map into $X( Y)$.
To correct this, let $U\subset \R^N$ be a tubular neighborhood
of $X( Y)$, and let $\Pi\colon U\to X( Y)$ be
a smooth retraction.
Define a $C^3$ map $Z\colon \overline M_R\to \overline M$ by
\begin{displaymath}
Z(x,\rho) = (X^{-1}\circ\Pi\circ \tw X(x,\rho),r(x,\rho)).
\end{displaymath}
The restriction of $Z$ to $\del M$ is the identity, and for some small
$\epsilon>0$, $Z$ is an embedding of $\overline M_\epsilon$ into $M$.

Let $p\in M_R$ be arbitrary, and let $q=Z(p)$.  Writing
$Z_*g = (Z^{-1})^*g$,
we conclude from \eqref{eq:dr2} that
\begin{align*}
|d\rho(q)|^2_{\rho^2 Z_*g}
&= \rho(q)^{-2} |d\rho(q)|^2_{Z_*g}\\
&= \rho(Z(p))^{-2} |d(\rho\circ Z)(p)|^2_{g}\\
&= r(p)^{-2} |dr(p)|^2_g\\
&= |dr(p)|^2_{r^2 g}
= 1+ O(\rho(p)^2) = 1 + O(\rho(q)^2).
\end{align*}
To check the mixed tangential/normal components of $\rho^2 Z_*g$,
let $(\theta^\alpha,\rho)$ be any background coordinates
on $\overline M_R$, and let $Z^\alpha = \theta^\alpha\circ Z$
denote the tangential component functions of $Z$ in these coordinates.
Observe that
$d\theta^\alpha(Z_*P) = d\theta^\alpha(X^{-1}_*\circ \Pi_* \circ \tw X_*P)$.
Because the component functions of $\Pi_*$ (as a map from $\R^N$
to itself) are uniformly
bounded, as are those of $X^{-1}_*$ in  background coordinates,
it follows that the tangential components $PZ^\alpha$ of $Z_*P$ in
background coordinates are  $O(\rho^2)$.
Thus
\begin{align*}
\< d\rho(q), d\theta^\alpha(q)\>_{\rho^2Z_*g}
&= \rho(q)^{-2}\<d\rho(q),d\theta^\alpha(q)\> _{Z_*g}\\
&= \rho(Z(p))^{-2}\<d(\rho\circ Z)(p),d(\theta^\alpha\circ Z)(p)\>_{g}\\
&= \<dr(p),dZ^\alpha(p)\>_{r^2 g}\\
&= P_pZ^\alpha = O(\rho(p)^2) = O(\rho(q))^2.
\end{align*}

We define our collar diffeomorphism by 
$G = Z^{-1}|_{M_{R_0}}$ for $R_0$
sufficiently small,
and let $\hat g = \rho^2 G^*g = \rho^2 Z_*g$.  By construction,
in any
background coordinates,
\begin{align}
|d\rho|^2_{\hat g} &= 1+O(\rho^2)\label{eq:drho2}\\
\<d\rho,d\theta^\alpha\>_{\hat g} &= O(\rho^2).
\end{align}
Inverting the coordinate matrix of $\hat g$, therefore,
we find that
\begin{displaymath}
\hat g = d\rho^2 +
\hat g_{\alpha\beta}(\theta,\rho)d\theta^\alpha d\theta^\beta + O(\rho^2)
\end{displaymath}
for some functions $\hat g_{\alpha\beta}$ that are $C^2$ up to $\del M$.
Moreover, because the restriction of $G$ to $\del M$ is the identity and 
$G^*\rho = \rho + O(\rho^2)$, 
$\hat g_{\alpha\beta}=\tw\gamma_{\alpha\beta} = \overline h_{\alpha\beta}$ at points
of $\del M$.

To conclude the proof, we will use the Einstein
equation to show that $\del_\rho\hat g_{\alpha\beta}=0$
along $\del M$, which implies $\hat g = \overline h + O(\rho^2)$
as desired.
In terms of $\hat g$, the Einstein equation for $G^* g$ translates to
\begin{displaymath}
-n \rho^{-2} \hat g_{jk} = \hat R_{jk} + (n-1)\rho^{-1} \rho_{;jk} + \rho^{-1}\rho_{;l}{}^l \hat g_{jk}
- n \rho^{-2} \rho_{;l} \rho_{;}{}^{l} \hat g_{jk},
\end{displaymath}
where the semicolons indicate covariant derivatives, all taken
with respect to $\hat g$ (cf.\/~\cite[p.~266]{Lee:spectrum}).
Multiplying by $\rho$, using \eqref{eq:drho2}, and evaluating at
$\rho=0$, we obtain $(n-1)\rho_{;jk} + \rho_{;l}{}^l\hat g_{jk}=0$
along $\del M$. Taking the trace with respect to $\hat g$, we find
that $\rho_{;l}{}^l=0$ and therefore $\rho_{;jk}=0$ along $\del
M$. Expanding this equation in terms of the Christoffel symbols of
$\hat g$ in background coordinates, we conclude that
$\del_\rho\hat
g_{\alpha\beta} = 0$ along $\del M$ as claimed.%
\end{proof}

\section{The harmonic map normalization}

In this section, we will show that $g$
can be modified by a collar diffeomorphism
so that it satisfies the elliptic equation
\eqref{eq:gauge-broken-Einstein} near the boundary.
We seek a collar diffeomorphism that is harmonic
from $(M_R,g)$ to $(M_R,h)$, where $h$ is the smoothly conformally
compact metric defined in the preceding section.
In order to find one, we will
parameterise the diffeomorphisms
near the identity by small vector fields using the Riemannian
exponential map of $h$.

For any small $R>0$, let $\del_R M_R =  Y\cross \{R\}$ 
denote
the ``inner boundary'' of $\overline M_R$, and let
$\Ccirc^{k,\lambda}_\delta(M_R;TM)$ denote the set of vector
fields in $C^{k,\lambda}_\delta(M_R;TM)$ that vanish on $\del_R
M_R$. If $v\in \Ccirc^{k,\lambda}_\delta(M_R;TM)$, define a map
$H_v\colon M_R\to M$ by
\begin{displaymath}
H_v(p) = \exp_p(v(p)),
\end{displaymath}
where $\exp$ denotes the Riemannian exponential map of $h$.
Since conformally compact metrics are complete 
at infinity \cite{mazzeo:hodge},
$H_v$ is well-defined as a map from $M_R$ into $M$ as long as both $R$ and
$v$ are sufficiently small.

Let us call a map $H\colon M_R\to M_R$
{\it admissible} if for each M\"obius chart, $H$ maps
$\overline{V_1(p_i)}$ into $V_2(p_i)$.  Because $h$ is
uniformly equivalent to the Euclidean metric in M\"obius coordinates,
for any admissible map $H$,
the Riemannian distance $d_h(p,H(p))$ is
uniformly bounded for $p\in M_R$.
This implies that $d_{\overline h}(p,H(p))\to 0$
uniformly as $p\to\del M$, which in turn implies that any
admissible map has a continuous extension to
$\del M$ that fixes $\del M$ pointwise.

\begin{lemma}\label{lemma:Hv-admissible}
If $\delta\ge 0$ and $v$ is sufficiently small in
$\Ccirc^{1,0}_{\delta}(M_R;TM)$, then $H_v$ is an admissible map
from $M_R$ to itself.
\end{lemma}

\begin{proof}
Because 
$d_h(p,H_v(p)) \le |v(p)|_h$ 
and $d_h$ is uniformly equivalent to Euclidean distance in
M\"obius coordinates, it follows that $H_v$ will be an admissible
map if $\|v\|_{0,0,\delta}$ is sufficiently small, provided that
$H_v$ maps $M_R$ to itself.  By examining the lengths of minimizing geodesics
to $\del_R M_R$, the reader can verify that this is the case provided that
$\|v\|_{1,0,\delta}$ is small enough that $|\Del v|_h\le \half$ on $M_R$.
\end{proof}

Let $\Sigma^2$ denote the bundle of symmetric covariant 2-tensors
over $M$. For any section $w$ of $\Sigma^2$, write $g_w = h +
w$. For any $0<\lambda<1$, define a map
\begin{displaymath}
\Theta\colon \Ccirc^{2,\lambda}_{1+\lambda}(M_R;TM)\cross C^{1,\lambda}_{1+\lambda}(M_R;\Sigma^2)\to
C^{0,\lambda}_{1+\lambda}(M_R;TM)\cross C^{1,\lambda}_{1+\lambda}(M_R;\Sigma^2)
\end{displaymath}
by
\begin{displaymath}
\Theta(v,w) = \left( (H_v)_*^{-1}(\Delta_{g_wh}H_v), w\right),
\end{displaymath}
where $\Delta_{g_wh}H_v$ denotes the harmonic map Laplacian of $H_v$,
viewed as a map from $(M_R,g_w)$ to $(M_R,h)$.

\begin{lemma}\label{lemma:Psi-C1}
The map $\Theta$ is well-defined and
of class $C^1$ in a neighborhood of $(0,0)$ in
$\Ccirc^{2,\lambda}_{1+\lambda}(M_R;TM)\cross C^{1,\lambda}_{1+\lambda}(M_R;\Sigma^2)$.
\end{lemma}

\begin{proof}
When $v$ and $w$ are understood, let us write $g=g_w$ and $H = H_v$ for brevity.
Recall that
$g$ is uniformly $C^{1,\lambda}$
equivalent to
the Euclidean metric in M\"obius coordinates;
a similar statement applies to $h$, but in that case we have
uniform $C^m$ equivalence for every $m$.

For this proof, we will denote M\"obius coordinates generically by
$x$ or $(x^j) = (x^1,\dots,x^{n+1})$, and the associated standard
fiber coordinates on $TM$ by $v$ or $(v^j) = (v^1,\dots,v^{n+1})$.
Letting $E^j(x,v)$ denote the
(smooth) component functions of the $h$-exponential map in
M\"obius coordinates, we see that $H$ has component functions
given by
\begin{displaymath}
H^j(x) = E^j(x,v(x)).%
\end{displaymath}
Because $E^j(x,0) = x^j$,
it follows from Lemma \ref{lemma:holder-facts}
that the functions $A^j(x) = H^j(x) - x^j$
satisfy the following uniform bound for sufficiently small $v\in \Ccirc^{2,\lambda}_{1+\lambda}(M_R;TM)$:%
\begin{displaymath}
\|A^j\|_{2,\lambda;\overline B_1} \le C
\|\Phi_i^*v\|_{2,\lambda;\overline B_1}\le C\rho(p_i)^{1+\lambda} \|v\|_{2,\lambda,1+\lambda}.
\end{displaymath}

Calculating in M\"obius coordinates, we have
\begin{equation}\label{eq:harmonic-map-in-coords}
(\Delta_{gh}H)^j =
g^{kl}\left(- \del_k \del _l H^j + \Gamma_{kl}^m \del_m H^j - (\Pi_{mq}^j\circ H) \del_k H^m \del_l H^q\right),
\end{equation}
where $ \Gamma_{kl}^m$ are the Christoffel symbols of $g$ and
$\Pi_{mq}^j$ are those of $h$. Note that we can write the
difference $\Gamma^m_{kl} - \Pi^m_{kl}$ as follows:
\begin{align*}
\Gamma^m_{kl} - \Pi^m_{kl}
&= \half g^{mj} \left( \del_k w_{lj} + \del_l w_{kj} - \del_j w_{kl}\right)\\
&\qquad + \half (g^{mj}-h^{mj}) \left( \del_k h_{lj} + \del_l h_{kj} - \del_j h_{kl}\right).
\end{align*}
By virtue of Lemma \ref{lemma:holder-facts} again, this time with
$f$ equal to the $mj$-component of the map taking an $(n+1)\cross
(n+1)$ matrix to its inverse, this last expression is in
$C^{0,\lambda}(\overline B_1)$, with $C^1$ dependence on $w$, and
satisfies an estimate of the form
\begin{displaymath}
\|\Gamma^m_{kl} - \Pi^m_{kl}\|_{0,\lambda;\overline B_1}\le C \rho(p_i) ^{1+\lambda} \|w\|_{1,\lambda,1+\lambda}.
\end{displaymath}

Now rewrite \eqref{eq:harmonic-map-in-coords} as follows:
\begin{align*}
(\Delta_{gh}H)^j
&=
g^{kl}\biggl(- \del_k \del _l A^j + (\Gamma_{kl}^m -\Pi_{kl}^m)\del_m H^j +
\Pi^m_{kl}\del_m A^j - \Pi^j_{ml} \del_k A^m \\
&\quad - \Pi^j_{mq} \del _k H^m \del_l A^q +
(\Pi_{mq}^j\circ \Id - \Pi_{mq}^j\circ H) \del_k H^m \del_l H^q\biggr).
\end{align*}
Another application of Lemma \ref{lemma:holder-facts} shows that
this expression is in $C^{0,\lambda}(\overline B_1)$, with $C^1$ dependence
on $v$ and $w$, and with $C^{0,\lambda}$ norm bounded by a multiple of
$\rho(p_i) ^{1+\lambda} \left(\|v\|_{2,\lambda,1+\lambda}
+ \|w\|_{1,\lambda,1+\lambda}\right)$.  Finally, since the pushforward map
$H_*^{-1}\colon T_{H(p)}M\to T_pM$ is represented by the inverse of the matrix
$\del H^j /\del x^k(p) = \delta^j_k + \del A^j /\del x^k(p)$, one last application
of Lemma \ref{lemma:holder-facts} shows that $\Theta$ is a $C^1$ map as claimed.
\end{proof}

\begin{lemma}\label{lemma:DPsi-isomorphism}
If $R$ is sufficiently small,
the differential $D\Theta_{(0,0)}$ is 
a Banach space
isomorphism from 
$\colon \Ccirc^{2,\lambda}_{1+\lambda}(M_R;TM)\cross C^{1,\lambda}_{1+\lambda}(M_R;\Sigma^2)$ to
$C^{0,\lambda}_{1+\lambda}(M_R;TM)\cross C^{1,\lambda}_{1+\lambda}(M_R;\Sigma^2)$.
\end{lemma}

\begin{proof}
At $(v,w)=(0,0)$, the differential of $\Theta$ can be computed
as follows:
\begin{align*}
D\Theta_{(0,0)}(v,w)
&= \bigg(\left.\pd{t}\right|_{t=0} \left(( H_{tv})_*^{-1}
(\Delta_{hh}H_{tv})\right)
+ \left.\pd{t}\right|_{t=0}(\Delta_{g_{tw}h}\Id), w\biggr)\\
&= ( Lv + Aw, w),
\end{align*}
where $L$ is the linearisation of the harmonic map Laplacian
$\Delta_{hh}$ about the identity map,
and $A$ is some first-order linear differential operator that is
bounded
from $C^{1,\lambda}_{1+\lambda}(M_R;\Sigma^2)$ to $C^{0,\lambda}_{1+\lambda}(M_R;TM)$.
Clearly this is invertible if and only if
$L\colon \Ccirc^{2,\lambda}_{1+\lambda}(M_R;TM)\to
C^{0,\lambda}_{1+\lambda}(M_R;TM)$
is invertible.

A computation shows that $L = \Del^*_h\Del_h - \Ric(h)$.
Because $\Ric(h)$ approaches $-n h$ at $\del M$, it is straightforward to
check that (in the terminology of \cite{Lee:fredholm})
the characteristic exponents of $L$ are
\begin{displaymath}
s = 0,\ n+2,\ \frac{n+2\pm\sqrt{n^2+8n}}{2}.
\end{displaymath}
It follows that
$L$ has indicial radius $R = (n+2)/2$,
and therefore by \cite[Theorem C and Section 7]{Lee:fredholm},
it is Fredholm as an operator from $C^{k+2,\lambda}_\delta(M;TM)$ to
$C^{k,\lambda}_\delta(M;TM)$ for all $k\ge 0$, $0<\lambda<1$, and
$-1<\delta<n+1$.
Moreover, \cite[Lemma 7.12]{Lee:fredholm} shows that for a $1$-form $u$
supported in $M_R$,
\begin{displaymath}
\left(u,\Del_h^*\Del_h u\right) \ge \left(\frac{n^2}{4} + 1 - \epsilon\right)\|u\|_{0,2}^2,
\end{displaymath}
where $\epsilon$ can be made as small as desired by taking $R$ small.
Since the operator $\Del_h^*\Del_h$ commutes with the
index-raising isomorphism between $1$-forms and vector fields, the same
result holds for vector fields.
It follows that $L \sim \Del_h^*\Del_h + n$
satisfies an a priori $L^2$ estimate of the form
\begin{displaymath}
\|v\|_{L^2} \le C\|Lv\|_{L^2} \quad \text{for all $v\in C^\infty_c(M_R;TM)$}
\end{displaymath}
when $R$ is sufficiently small.
Then the same argument as in the proof of Theorem C of \cite{Lee:fredholm}
implies that
$L\colon \Ccirc^{k+2,\lambda}_\delta(M_R;TM)\to
C^{k,\lambda}_\delta(M_R;TM)$ is an isomorphism
for $-1<\delta<n+1$; the only modification that
needs to be made is to handle the
Dirichlet boundary condition on the
inner boundary $\del_R M_R$, but as $L$ is uniformly elliptic there, the
required estimates follow easily from the standard theory of elliptic
boundary value problems.
\end{proof}

Now suppose $g$ satisfies the hypotheses of Theorem \ref{thm:main}, and
let $w = G^*g-h$, where $G$ is given by Lemma \ref{lemma:bdry-adjusting-diffeo}.
Let $\psi\colon \R\to
[0,1]$ be a smooth cutoff function such that $\psi(t) \equiv 1$
for $t\le \half$ and $\psi(t)\equiv 0$ for $t\ge 1$. For any small
$s>0$, define $\psi_s\in C^\infty(\overline M_R)$ by
\begin{displaymath}
\psi_s(p) = \psi\left(\frac{\rho(p)}{s}\right),
\end{displaymath}
Then we define $w_s = \psi_s w$.
Observe that $g_{w_s} = G^*g$ on
the subset $M_{s/2}$ where $\psi_s\equiv 1$.

\begin{lemma}\label{lemma:w_s->0}
For any fixed small $R>0$ and any $0<\lambda<1$,
$w_s\to 0$ in $C^{1,\lambda}_{1+\lambda}(M_R;\Sigma^2)$ as $s\to 0$.
\end{lemma}

\begin{proof}
Because $\psi_s$ is
uniformly bounded in $C^{1,\lambda}(M_R)$ and supported in $M_s$, 
the lemma follows from the fact 
that, in M\"obius coordinates, the component
functions of $\Phi^*_iw$  and their first and second derivatives are bounded by
a constant multiple of $\rho(p_i)^2$.
\end{proof}

\begin{theorem}\label{thm:harmonic-map-existence}
With $g$ and $h$ as above,
for any $0<\lambda<1$,
there exists
a $C^{2,\lambda}$ collar
diffeomorphism $\Psi \colon \overline M_R\to \overline M$ such that
$\Psi ^*g-h\in C^{1,0}_{1+\lambda}(M_R;\Sigma^2)$
and $\tw g = \Psi ^*g$
satisfies
\eqref{eq:gauge-broken-Einstein} on $M_{R_0}$ for some $0<R_0<R$.
\end{theorem}

\begin{proof}
Let $\Theta$ and $w_s$ be defined as above.
It follows from Lemma
\ref{lemma:DPsi-isomorphism} and the inverse function theorem
that $\Theta$ is a bijection from
a neighborhood of $(0,0)$ in $\Ccirc^{2,\lambda}_{1+\lambda}(M_R;TM)
\cross C^{1,\lambda}_{1+\lambda}(M_R;\Sigma^2)$
to a neighborhood of $(0,0)$ in
$C^{0,\lambda}_{1+\lambda}(M_R;TM)\cross C^{1,\lambda}_{1+\lambda}(M_R;\Sigma^2)$.
By Lemma \ref{lemma:w_s->0}, therefore, we can choose $s$ small enough that
$(0,w_s) = \Theta(v,w)$ for some $(v,w)\in \Ccirc^{2,\lambda}_{1+\lambda}(M_R;TM)
\cross C^{1,\lambda}_{1+\lambda}(M_R;\Sigma^2)$.
This is equivalent to the assertion that $w=w_s$ and
$H_v$ is harmonic from $(M_R,g_{w_s})$ to $(M_R,h)$.
Because $g_{w_s} = G^*g$ on $M_{s/2}$,
$H_v$ is also harmonic from $(M_{s/2},G^*g)$ to $(M_{s/2},h)$.

Because the
component functions $H^j(x)$ in M\"obius coordinates differ from
$x^j$ by functions $A^j(x)$ that can be made as small as desired
in $C^{2,\lambda}(\overline B_2)$ (by taking $s$ sufficiently
small), it follows that $H_v$ is a diffeomorphism from $M_R$ to
itself, and by Lemma \ref{lemma:Hv-admissible} it is an admissible
map and therefore extends to a homeomorphism of $\overline M_R$
fixing $\del M$ pointwise, i.e., a collar diffeomorphism.

Define $\Psi = G\circ H_v^{-1}\colon M_R\to M$, and
let $\tw g = \Psi ^* g$.  By
the diffeomorphism invariance of the Einstein equation, $\tw g$ is Einstein; and by
the diffeomorphism invariance of the harmonic
map equation, the identity map is harmonic
on $M_{s/2}$ from
$\tw g = \left(H_v^{-1}\right)^*(G^*g)$ to $h$.  Thus
$\tw g$
satisfies
\eqref{eq:gauge-broken-Einstein} on $M_{s/2}$.

Next we will show that $\Psi ^*g-h\in
C^{1,0}_{1+\lambda}(M_R;\Sigma^2)$. Since $G^*g - h = w\in
C^{1,\lambda}_{1+\lambda}(M_R;\Sigma^2)$ by the proof of Lemma
\ref{lemma:w_s->0}, it suffices to show that $(H_v^{-1})^* G^*g -
G^*g\in C^{1,\lambda}_{1+\lambda}(M_R;\Sigma^2)$. Let us
abbreviate $G^*g$ by $\hat g$ and $H_v^{-1}$ by $K\colon M_R\to
M_R$.  Again, by taking $s$ (and thus also $A^j$) small enough, we
can ensure that $H_v(V_2(p_i))$ contains $\overline V_1(p_i)$ for
each M\"obius chart, and therefore $K$ is an admissible map.  We
will write the component functions of $H_v$ as $H^j(x) = x^j +
A^j(x)$ as above, with
\begin{displaymath}
\|A^j\|_{2,\lambda;\overline B_2} \le C\rho(p_i)^{1+\lambda}.
\end{displaymath}

If $s$ is chosen small enough, then
the functions $A^j$ will be uniformly small in $C^{2,\lambda}(\overline B_2)$, so
the Jacobian matrix $\del_k H^j(x)$ will be uniformly invertible on $\overline B_2$,
independent of the choice of M\"obius chart.
Let us write $(I^j_k(x))$ for the components of the inverse matrix
of $(\del_k H^j(x))$.
Because matrix inversion is
continuous in the $C^{1,\lambda}$ norm by
Lemma \ref{lemma:holder-facts},
the functions $I^j_k$ are uniformly bounded in $C^{1,\lambda}(\overline B_2)$.

The chain rule shows that
\begin{equation}\label{eq:chain-rule-inverse}
\del_k K^j(x) = I^j_k(K(x)),
\end{equation}
which is continuous and uniformly bounded on $\overline B_1$.
But this implies that $K^j$ is uniformly $C^1$, and  using
\eqref{eq:chain-rule-inverse} we conclude successively that
$\del_k K^j(x)$ is uniformly $C^1$, $K^j$ is uniformly $C^2$,
$\del_k K^j(x)$ is uniformly $C^{1,\lambda}$ (by Lemma \ref{lemma:holder-facts},
and thus $K^j$ is
uniformly $C^{2,\lambda}$ on $\overline B_1$.

Let us write $B^j(x) = K^j(x) - x^j$.
The fact that
$H_v\circ K=\Id$ translates in M\"obius coordinates to
\begin{displaymath}
K^j(x) + A^j(K(x)) = x^j,
\end{displaymath}
so $B^j(x) = -A^j(K(x))$,
and we conclude that
\begin{equation}\label{eq:B_j-bound}
\|B^j\|_{2,\lambda;\overline B_1}
= \|A^j\circ K\|_{2,\lambda;\overline B_1}
\le C \|A^j\|_{2,\lambda;\overline B_2} \le C'\rho(p_i)^{1+\lambda}.
\end{equation}

Now
to show that $K^*\hat g - \hat g \in C^{1,0}_{1+\lambda}(M_R;\Sigma^2)$,
we just compute in M\"obius
coordinates:
\begin{align*}
K^*\hat g - \hat g
&= \hat g_{jk}(K(x)) (dx^j+ dB^j)(dx^k+dB^k) - \hat g_{jk}(x)dx^j\,dx^k\\
&= \Bigl(\hat g_{jk}(K(x)) - \hat g_{jk}(x)\Bigr) dx^j\, dx^k \\
&\qquad + 2 \hat g_{jk}(K(x))\frac{\del B^k}{\del x^q}dx^j\,dx^q
+ \hat g_{jk}(K(x)) \frac{\del B^j}{\del x^m}\frac{\del B^k}{\del x^q}dx^m\,dx^q.
\end{align*}
Because $\rho^2\hat g\in C^2(\overline M)$, the component functions
$\hat g_{jk}$ are uniformly bounded in $C^2(\overline B_2)$ by the
properties of M\"obius coordinates.  The same is true of $\hat g_{jk}\circ K$
by composition.  The last two terms above thus have $C^1$ norms uniformly
bounded by a constant multiple of $\rho(p_i)^{1+\lambda}$.
Differentiating the
first term, we obtain%
\begin{align*}
\del_l \left( \hat g_{jk}\circ K - \hat g_{jk}\right)
&= \left(\del_m \hat g_{jk} \circ K\right) \frac{\del K^m}{\del x^l}
- \del_l \hat g_{jk} \\
&= \left(\del_m \hat g_{jk} \circ K\right) \frac{\del B^m}{\del x^l}
+ \left(\del_l \hat g_{jk} \circ K -\del_l \hat g_{jk}\right).
\end{align*}
The first term is uniformly bounded by a multiple of $\rho(p_i)^{1+\lambda}$
thanks to \eqref{eq:B_j-bound}.  The same is true of the
second term by a simple application of the mean value theorem and
the fact that $K^j(x) - x^j = B^j(x)$.
This completes the proof that
$\Psi ^*g-h\in C^{1,0}_{1+\lambda}(M_R;\Sigma^2)$.

It remains only to show that
$\Psi$ has a $C^{2,\lambda}$ extension to $\overline M_R$.
Because $G$ is $C^3$ on $\overline M_R$ by construction, it suffices
to consider $H_v^{-1}= K$.
Choose a fixed $p_i\in M_R$ 
and corresponding
M\"obius chart $\Phi_i$.
We can write the map $K$ either in 
M\"obius coordinates, with coordinate functions denoted by 
$(K^j)$ as above: 
\begin{displaymath}
K(x^j) = \left(K^1(x^j),\dots,K^{n+1}(x^j)\right),
\end{displaymath}
or in background coordinates $(\theta^j) = (\theta^1,\dots,\theta^n,\rho)$,
with coordinate functions that we will denote by $(\Kbar^j)$:
\begin{displaymath} 
K(\theta^j) = \left(\Kbar^1(\theta^j),\dots,\Kbar^{n+1}(\theta^j)\right).
\end{displaymath}
The two coordinate representations are related by
\begin{displaymath}
\Kbar^m(\theta^j) = c^m + \rho(p_i)K^m\left(\frac{\theta^j - c^j}{\rho(p_i)}\right),
\end{displaymath}
where $c^j$ are constants defined by 
$(c^1,\dots,c^{n+1}) = (\theta^1(p_i),\dots,\theta^n(p_i),0)$.
Using once again the fact that $K^j(x) = x^j+ B^j(x)$ with $B^j$ satisfying
\eqref{eq:B_j-bound}, we compute
\begin{displaymath}
\frac{\del \Kbar^m}{\del \theta^k}
= \delta^m_k + \frac{\del B^m}{\del x^k}\left( \frac{\theta^j - c^j}{\rho(p_i)}\right),
\end{displaymath}
so both $\Kbar^m$ and $\del \Kbar^m/\del x^k$ are uniformly bounded.
Differentiating once more, we find
\begin{displaymath}
\frac{\del^2 \Kbar^m}{\del \theta^k\del \theta^l}
= \rho(p_i)^{-1}\frac{\del^2 B^m}{\del x^k\del x^l}\left( \frac{\theta^j - c^j}{\rho(p_i)}\right),
\end{displaymath}
and therefore,
\begin{align*}
\left| \frac{\del^2 \Kbar^m}{\del \theta^k\del \theta^l}(\theta^j) \right.&-
\left.\frac{\del^2 \Kbar^m}{\del \theta^k\del \theta^l}(\tw\theta^j) \right|\\
&= \rho(p_i)^{-1}\left|
\frac{\del^2 B^m}{\del x^k\del x^l}\left( \frac{\theta^j - c^j}{\rho(p_i)}\right)
- \frac{\del^2 B^m}{\del x^k\del x^l}\left( \frac{\tw\theta^j - c^j}{\rho(p_i)}\right)\right|\\
&\le \rho(p_i)^{-1}
\|B^m\|_{2,\lambda;\overline B_2}
\left|\frac{\theta^j - \tw\theta^j}{\rho(p_i)}\right|^\lambda\\
&\le C
\left|\theta^j - \tw\theta^j\right|^\lambda,
\end{align*}
which shows that $\Kbar^m$ is uniformly $C^{2,\lambda}$ up to the boundary
as claimed.
\end{proof}

\section{Polyhomogeneity}\label{section:polyhomogeneity}

Let $U_0\subset \R^n$ be an open set, and let 
$U= U_0\cross(0,\epsilon)
\subset \H^{n+1}$.  For any $\delta\in\R$, 
we denote by $\scr C^\delta$ the space
of functions $f\in C^\infty(U)$ that satisfy,
on any subset $K\cross(0,\epsilon_0)$ with
$K\subset U_0$ compact and $0<\epsilon_0<\epsilon$, 
estimates of the following form
for all integers $r\ge 0$ and all multi-indices $\alpha$:
\begin{displaymath}
\left| (y\del_y)^r \del_ x^\alpha f(x,y)\right| \le C_{r,\alpha} y^\delta.
\end{displaymath}
(We use the multi-index notations 
$\alpha = (\alpha_1,\dots,\alpha_{n})$ and
$\del_ x^\alpha = 
(\del_{ x^1})^{\alpha_1}\dots (\del_{ x^n})^{\alpha_n}$.)

A smooth function $f\colon U\to \R$ is said to be {\it
polyhomogeneous} (cf.\/~\cite{mazzeo,AndChDiss}) if there
exists a sequence of real numbers $s_i \to +\infty$, a sequence of
nonnegative integers $\{q_i\}$, and  functions $f_{ij}\in
C^\infty(U_0)$  such that
\begin{equation}\label{eq:def-polyhomo}
f(x,y) \sim \sum_{i=1}^\infty \sum_{j=0} ^{q_i} y^{s_i} (\log y)^j f_{ij}(x)
\end{equation}
in the sense that
for any $\delta>0$, there exists a positive
integer $N$ such that 
\begin{displaymath}
f(x,y) -
\sum_{i=1}^{N}\sum_{j=0}^{q_i} y^{s_i} (\log y)^j f_{ij}(x)
\in \scr C^\delta.
\end{displaymath}
A function or tensor field on $M_R$ is said to be polyhomogeneous if its coordinate
representation in every background chart is polyhomogeneous.
(Note 
that the definition of polyhomogeneity is phrased somewhat
differently in ~\cite{AndChDiss}, but it is easy to verify that the two
definitions are equivalent.  Note also that there is a misprint in the 
first displayed inequality of \cite[Section 3.1.5]{AndChDiss}: $\del_y^\alpha$
in that inequality
should be replaced by $\del_v^\alpha$.)

In this section, we will
apply the theory of \cite{AndChDiss} to conclude that solutions to
\eqref{eq:gauge-broken-Einstein} are polyhomogeneous.
A key step in the proof will be a regularity result for the linearised operator
$D_1 Q_{(h,h)}=\half(\Delta_L + 2n)$ from \cite{Lee:fredholm}.
Following \cite{AndChDiss},
we say that an interval $(\delta_-,\delta_+)\subset\R$ is a
{\it \(weak\) regularity interval}
for a second-order linear operator
$P$ on the spaces $C^{k,\lambda}_\delta(M_R;\Sigma^2)$
if whenever $u$ is a locally $C^2$ section of $\Sigma^2$ such that
$u\in C^{0,0}_{\delta_0}(M_R;\Sigma^2)$
and $Pu\in C^{0,\lambda}_{\delta}(M_R;\Sigma^2)$ with $\lambda\in (0,1)$
and $\delta_-<\delta_0<\delta<\delta_+$, it follows that
$u\in C^{2,\lambda}_{\delta}(M_R;\Sigma^2)$.
(We caution the reader that the notations for weighted H\"older spaces
used in \cite{AndChDiss} are different from those of \cite{Lee:fredholm}
that we are using here, in terms of index positions and
normalization of weights.  The space of sections of $\Sigma^2$
that is denoted by $C^\delta_{k+\lambda}(M)$
in \cite{AndChDiss}
is equal to the space that we would call $C^{k,\lambda}_{\delta+2}(M_R;\Sigma^2)$.
The difference between the weight factors $\delta$ and $\delta+2$
arises because
we measure the size of the component functions in M\"obius coordinates,
while \cite{AndChDiss} measures them in background coordinates.
To avoid confusion, we will use the notation $AC^\delta_{k+\lambda}(M)$ for the space
denoted by $C^\delta_{k+\lambda}(M)$ in \cite{AndChDiss}, so that
$AC^\delta_{k+\lambda}(M)=C^{k,\lambda}_{\delta+2}(M_R;\Sigma^2)$.
The condition that we have defined here would be expressed in \cite{AndChDiss}
by saying that $(\delta_--2,\delta_+-2)$ is a regularity interval
for $P$ on the spaces $AC^\delta_{k+\lambda}(M)$.)

\begin{theorem}
With $\tw g = \Psi ^*g$ as in the preceding section,
$\tw g$ is polyhomogeneous.
\end{theorem}

\begin{proof}
For any small symmetric $2$-tensor $\phi$ on $M_R$, define
\begin{displaymath}
F[\phi] := \rho^2 Q(h+\rho^{-2}\phi, h),
\end{displaymath}
with $Q$ as in \eqref{eq:gauge-broken-Einstein}.
Then $\phi = \rho^2(\tw g - h)$ satisfies $F[\phi] = 0$.
We wish to apply \cite[Theorem 5.1.1]{AndChDiss} to $F$, and thereby conclude that
$\phi$ is polyhomogeneous.
To do so, we will consider $\phi_0\equiv 0$ as an approximate
solution to $F[\phi]=0$, and check that
$F$, $\phi$, and $\phi_0$
satisfy each of the hypotheses of Theorem 5.1.1.
(Actually, we will be using the modified version of Theorem 5.1.1 described 
in Remark (ii) following the statement of that theorem---we have only a 
{\it weak} regularity interval, but we will verify that our constants 
all satisfy appropriately strengthened inequalities so that this is sufficient.)

\begin{enumerate}\romanletters
\item
{\it $F$
is a geometric operator in the sense of
\cite{AndChDiss}}: This follows easily from the fact that
$F$ is an invariant operator on tensor fields.
(Note that the notion of geometricity defined in \cite{AndChDiss} is
considerably weaker than that of \cite{Lee:fredholm}; in particular,
the definition in \cite{AndChDiss} does not require the coordinate
expression of $F[\phi]$ to depend only on the coefficients of
a single metric and the covariant derivatives of its curvature.)
\item
{\it $F$ is quasilinear, and $F[\phi]$ can be written in background
coordinates as a smooth function of
$\left(\rho, \theta^\alpha, \phi_{ij}, \rho\del_{k}\phi_{ij},
\rho^2\del_{k}\del_{l}\phi_{ij}\right)$}:
This is easily seen by expanding $Q(\tw g,h)$ in background coordinates.
\item
{\it $F[\phi_0]$ is a smooth tensor field on $\overline M$ and an
element of $AC^{\delta_0}_0(M) = C^{0,0}_{\delta_0+2}(M;\Sigma^2)$
for some $\delta_0>1$}: Observe that $F[\phi_0] = F[0] = \rho^2
Q(h,h) = \rho^2(\Ric(h) + n h)$. Using the formula for the
transformation of the Ricci tensor under a conformal change of
metric (cf.\/~\cite[p.\ 59]{Besse}) this can be written in
background coordinates as
\begin{align*}
F[\phi_0]&= \rho^2\left(R_{ij} + n h_{ij}\right) \\
&= \rho^2\overline R_{ij} + (n-1) \rho \rho_{;ij}
+ \rho \overline h^{kl}\rho_{;kl}\overline h_{ij} - n\overline h^{kl}\rho_{;k}\rho_{;l} \overline h_{ij}
+ n \overline h_{ij},
\end{align*}
where the curvature and covariant derivatives are computed with respect to the
smooth metric $\overline h = \rho^2 h$.
Because of the way we constructed $\overline h$,
$|d\rho|_{\overline h}^2$ is identically
equal to $1$ and $\rho_{;ij}$ is identically zero, so
all terms after the first one cancel, showing that $F[\phi_0] = \rho^2 \overline R_{ij}$, which
is smooth on $\overline M$ and $O(\rho^2)$ in background coordinates, and thus
an element of $AC^{2}_0(M)$.
\item
{\it $\phi - \phi_0\in AC^\delta_1(M_R)= C^{1,0}_{\delta+2}(M_R;\Sigma^2)$ for some $\delta>1$
and $R>0$}:
Since $\phi-\phi_0 = \phi = \rho^2 (\tw g - h)$, this is
equivalent to the assertion that $\tw g - h\in C^{1,0}_{\delta}(M_R;\Sigma^2)$
for some $\delta>1$,
which is guaranteed by Theorem \ref{thm:harmonic-map-existence}.
\item
{\it The linearised operator $F'[\phi_0]$
is a geometric elliptic
operator satisfying conditions {\rm (4.2.1)--(4.2.4)} of \cite{AndChDiss}}:
Actually, this is not true as stated, but something just as good is
true.
Note that $F'[\phi_0] = \rho^{2}\circ D_1 Q_{(h,h)}\circ\rho^{-2}
= \rho^{2}\circ(\Delta_L+2n)\circ\rho^{-2}$, which is certainly a
geometric elliptic operator.

Define subbundles of $\Sigma^2$ as follows:
\begin{align*}
V_0 &= \vectorspan(g);\\
V_1 &= \left\{q\in \Sigma^2: \Tr_g q = 0,\ q(\grad\rho,\centerdot) =0\right\};\\
V_2 &= \vectorspan\left( g - \frac{n+1}{|d\rho|_g^2} d\rho\otimes d\rho\right);\\
V_3 &= \left\{d\rho\otimes\omega + \omega\otimes d\rho: \omega\in T^*M,\ \left<\omega,d\rho\right>= 0\right\}.
\end{align*}
It is easy to check that $\Sigma^2$ admits an orthogonal decomposition $\Sigma^2=V_0\oplus V_1 \oplus
V_2\oplus V_3$.  For $i=0,\dots,3$, let $\pi_i\colon \Sigma^2\to V_i$ denote the
orthogonal projection.
The arguments of \cite[pp.\ 199--202]{GL} show that $\Delta_L+2n$ can be written
in the form
\begin{displaymath}
\Delta_L + 2n = \sum_{i=0}^3 P_i + \tw P,
\end{displaymath}
where in background coordinates
$\tw P$ is of the form \cite[(4.2.3)]{AndChDiss} and each $P_i$ is an operator on
sections of $V_i$ that
can be written
\begin{displaymath}
P_i  = -\left( \rho^2 \del_\rho^2 + (5-n) \rho\del_\rho + B_i\right) \otimes \pi_i,
\end{displaymath}
with $B_0 = B_2 = 4-2n$, $B_1 = 4-4n$, and $B_3 = 3-3n$.
Thus
\begin{displaymath}
F'[\phi_0] = \sum_{i=0}^3 L_i + \tw L,
\end{displaymath}
where $\tw L= \rho^2\circ \tw P \circ \rho^{-2}$ is again of the form \cite[(4.2.3)]{AndChDiss}, and
\begin{displaymath}
L_i = \rho^{2} \circ P_i \circ \rho^{-2} =
-\left( \rho^2 \del_\rho^2 + (1-n) \rho\del_\rho + b_i\right) \otimes \pi_i,
\end{displaymath}
with $b_i = B_i-4< 0$.
All of the arguments in \cite[Sections 4 and 5]{AndChDiss}
go through with only trivial changes if
the ordinary differential operator denoted there by $L_{ab}$
is replaced by the block-diagonal operator $L_0\oplus L_1\oplus L_2\oplus L_3$.
\item
{\it The interval $(0,n)$ is a regularity interval for the
operator $F'[\phi_0]$ on the spaces $AC^\delta_{k+\lambda}(M_R)$}:
It is an immediate consequence of Lemma 6.4(b) of
\cite{Lee:fredholm} that $(0,n)$ is a regularity interval for
$\Delta_L+2n$ on the spaces $C^{k,\lambda}_\delta(M;\Sigma_R)$,
and it follows immediately from this that $(-2,n-2)$ is a
regularity interval for $F'[\phi_0] = \rho^{2}\circ
(\Delta_L+2n)\circ\rho^{-2}$ on the same spaces. In the
terminology of \cite{AndChDiss}, this means that $(0,n)$ is a
regularity interval on $AC^\delta_{k+\lambda}(M_R)$.
\end{enumerate}

From \cite[Theorem 5.1.1]{AndChDiss}, therefore, we conclude that
$\phi$ (and hence also $\tw g$) is polyhomogeneous.
\end{proof}

\section{The asymptotic expansion}

To obtain the asymptotic expansion announced in 
Theorem \ref{thm:main}, we need to subject $g$ to one more
collar diffeomorphism.  First, a preliminary lemma.%

\begin{lemma}\label{lemma:existence-of-r}
Suppose $\tw g$ is a Riemannian metric on $M_R$
that is polyhomogeneous and conformally compact of class $C^{1,\lambda}$ for some
$0<\lambda<1$, and satisfies $|d\rho|_{\rho^2\tw g}\to 1$ at $\del M$.
Then there exists a polyhomogeneous $C^{1,\lambda}$
defining function $r$ such that 
$|dr|_{r^2 \tw g}\equiv 1$ in a neighborhood of $\del M$
and $r/\rho\to 1$ at $\del M$.
\end{lemma}

\begin{proof}
Writing $\overline g =
\rho^2 \tw g$ and $r = \rho e^u$, we see that 
the conclusion is equivalent to
\begin{displaymath}
|d\rho|^2_{\overline g} + 2\rho\< d\rho,du\>_{\overline g} + \rho^2 |du|^2_{\overline g} = 1,
\qquad u|_{\del M} = 0.
\end{displaymath}
It is shown in \cite[Lemma 5.1]{Lee:spectrum} that this has a solution
$u\in C^{2,\lambda}(\overline M_R;\R)$ if $\overline g$ is $C^{3,\lambda}$ up to
the boundary, by reducing it to finding the flow of the
Hamiltonian vector field $X_F$, where $F\colon T^*\overline M_R\to\R$ is the function
defined by
\begin{displaymath}
F(\theta,\xi) = 2\<d\rho,\xi\>_{\overline g} + \rho |\xi|_{\overline g}^2 -
\frac{1-|d\rho|^2_{\overline g}}{\rho}.
\end{displaymath}
(Here $(\theta,\xi) =
(\theta^1,\dots,\theta^{n+1},\xi_1,\dots,\xi_{n+1})$ are standard
coordinates on $T^*\overline M_R$ associated with background coordinates.) 
In the present situation, $F$
is only of class $C^{0,\lambda}$ on $T^*\overline M_R$;  but because
it is polyhomogeneous, each of the following quantities is also 
$C^{0,\lambda}$ in each background chart:
\begin{gather*}
\frac{\del{F}}{\del \xi_j},\
\frac{\del{F}}{\del \theta^\alpha},\
\rho\frac{\del{F}}{\del \rho},\\
\frac{\del^2{F}}{\del \xi_j\del \xi_k},\
\frac{\del^2{F}}{\del \theta^\alpha\del\theta^\beta},\
\frac{\del^2{F}}{\del \xi_j\del \theta^\alpha},\
\rho\frac{\del^2{F}}{\del\theta^\alpha\del \rho},\
\rho\frac{\del^2{F}}{\del\xi_j\del \rho},\
\rho^2\frac{\del^2{F}}{\del\rho^2}.
\end{gather*} Since the normal component of
$X_F$ satisfies
\begin{displaymath}
d\rho(X_F)|_{\del T^*\overline M_R} = \left.\frac{\del F}{\del
\xi_{n+1}}\right|_{\del T^*\overline M_R} = 2\overline
g^{n+1,n+1}|_{\del   T^*\overline M_R}= 2,
\end{displaymath}
it follows that $V=\half X_F$ satisfies the hypotheses of Lemma
\ref{lemma:polyhomogeneous-ode} below, so the flow-out by $X_F$
from the boundary
of $T^*\overline M_R$ exists and is polyhomogeneous and of class $C^{0,\lambda}$.
The rest of the
argument in \cite[Lemma 5.1]{Lee:spectrum} then goes through to
prove the existence of $r$ as claimed. Moreover, the solution $r$
so obtained is itself polyhomogeneous; matching lowest-order terms
in the expansion of $r$ with those in the expansion of $\overline
g$, we find that $r\in C^{1,\lambda}(\overline M_R)$ as claimed.
\end{proof}

It is clear 
from the proof of the preceding lemma that
near any boundary point, one can choose new coordinates in 
which $r$ is one of the coordinate functions 
and the metric is polyhomogeneous.  In fact, we can do
much better, as the next lemma shows.

\begin{lemma}\label{lemma:phg-Fermi-coords}
If $\tw g$ satisfies the hypotheses of the preceding lemma, then for
$R$ sufficiently small,
there exists a polyhomogeneous $C^{1,\lambda}$ collar diffeomorphism
$\Gamma\colon \overline M_R\to \overline M$ 
such that $\Gamma^*\tw g$ has the form \eqref{eq:metric-in-Fermi-coords}.
\end{lemma}

\begin{proof}
Let $r$ be the polyhomogeneous $C^{1,\lambda}$ defining function
given by Lemma \ref{lemma:existence-of-r}.
Let $q\in\del M$, and let
$(\theta^\alpha,\rho)$ be a fixed choice of background coordinates
on a neighborhood $\Omega$ of $q$ in $\overline M_R$.
As in the proof of Lemma
\ref{lemma:bdry-adjusting-diffeo},
let $P$ be the $r^2 \tw g$-gradient of $r$, which is a polyhomogeneous $C^{0,\lambda}$
vector field on $\overline M_R$ whose normal component satisfies
$d\rho(P) = 1$ along $\del M$.
Using Lemma \ref{lemma:polyhomogeneous-ode} again, in each background coordinate
chart $\Omega$ we obtain a uniquely determined $C^{0,\lambda}$
polyhomogeneous
flow $(x,t)\mapsto \gamma_x(t)$, where for each $x\in \Omega\cap \del M$, 
$\gamma_x\colon [0,\epsilon]\to \overline M_R$ is the integral curve
of $P$ starting at $x$.  
Comparing lowest-order terms in the expansions of $\gamma_x(t)$
and $P$, we see that the flow is in fact $C^{1,\lambda}$ up to the
boundary.
The various maps thus obtained in different coordinate charts all agree
where they overlap, so they patch together to define a global map
$\Gamma\colon  \del M\cross [0,\epsilon]\to \overline M_R$. 
The inverse function theorem shows that 
$\Gamma$ is a diffeomorphism in a neighborhood of $\del M\cross\{0\}$.
Identifying $\del M\cross [0,\epsilon]$ with $M_\epsilon$, we can
view $\Gamma$ as a collar diffeomorphism.
It is easy to check that $\Gamma^*\tw g$ 
has the form \eqref{eq:metric-in-Fermi-coords} with 
$G(\rho)$ polyhomogeneous.  
\end{proof}

Here is the ODE lemma used in the proofs of Lemmas \ref{lemma:existence-of-r}
and \ref{lemma:phg-Fermi-coords}.
It is adapted
from Proposition B.1 in \cite{ChMS}.
In this lemma, $\H^{m+1}$ denotes the upper half-space in $\R^{m+1}$, with coordinates
$(x^1,\dots,x^m,y)$.
In our application of this lemma 
in the proof of Lemma \ref{lemma:existence-of-r}, the
$x^i$-coordinates correspond to $(\theta^1,\dots,\theta^n,\xi_1,\dots,\xi_{n+1})$,
and $y$ corresponds to $\rho = \theta^{n+1}$.

\begin{lemma}\label{lemma:polyhomogeneous-ode}
Let $U_0$ be an open subset of $\R^m$, and let 
$V$ be a $C^1$ vector field on $U = U_0\cross(0,\epsilon)\subset {\H^{m+1}}$
of the form
\begin{displaymath}
V = A^i(x,y)\pd{x^i} + (1 + B(x,y))\pd{y}.
\end{displaymath}
Suppose that $A^i$ and $B$ satisfy the following estimates for some constants
$C_0>0$ and $0<\lambda<1$:
\begin{equation}\label{eq:AB-estimates}
\begin{aligned}
|B|,\ |\del_{x^j} B|, |y\del_y B|  &\le C_0 y^\lambda,\\
|A^i|,\ |\del_{x^j}A^i|,\ |y\del_y A^i|&\le C_0y^{\lambda-1}.
\end{aligned}
\end{equation}
\begin{enumerate}\letters
\item
If $K$ is any compact subset of $U_0$,
there exists $\epsilon>0$ such that, for each $x_0\in K$, there is a
unique continuous solution $\gamma = \gamma_{x_0}$ on $[0,\epsilon]$
to the initial value problem
\begin{equation}\label{eq:gamma-ode}
\begin{aligned}
\gamma'(t) &= V(\gamma(t)),\\
\gamma(0) &= (x_0,0).
\end{aligned}
\end{equation}
\item
If the coefficient functions $A^i$ and $B$ are polyhomogeneous,
then the map $(x,t)\mapsto \gamma_{x}(t)$ is polyhomogeneous 
and $C^{0,\lambda}$ up
to the boundary.  
\end{enumerate}
\end{lemma}

\begin{proof}
Let $x_0\in K$ be arbitrary.  Choose constants $C_1>0$ and
$0<\alpha<\lambda$, and let
$\epsilon$ be a positive constant to be determined later.
Let $\scr X$ denote
the set of continuous maps $\gamma\colon[0,\epsilon]\to \R^{m+1}$
of the form
\begin{equation}\label{eq:gamma-a-b}
\gamma(t) = (x_0 + a(t), t + b(t))
\end{equation}
satisfying
\begin{align}
|a^i(t)| &\le C_1 t^\alpha,\label{eq:a-estimate}\\
|b(t)| &\le C_1 t^{\alpha+1}.\label{eq:b-estimate}
\end{align}
If $\epsilon$ is sufficiently small, all such maps take their values in $U$
for $t\in(0,\epsilon]$, and
$\scr X$ is a complete metric space when endowed with the metric
\begin{displaymath}
d(\gamma,\tw\gamma) =
\sup_{t\in[0,\epsilon]} \left( t^{-\alpha} |a(t)-\tw a(t)| + t^{-\alpha-1} |b(t)-\tw b(t)|\right).
\end{displaymath}
Note that for any $\gamma\in\scr X$, we have
$|y(\gamma(t)) - t| = |b(t)|\le C_1t^{1+\alpha}\le C_1\epsilon^\alpha t$, so
choosing $\epsilon$ small enough that $C_1\epsilon^\alpha<\half$ implies that
\begin{equation}\label{eq:b-gt-half-t}
\half t \le y(\gamma(t))\le \tfrac32 t.
\end{equation}

Define a map $T\colon \scr X\to \scr X$ by
\begin{displaymath}
T\gamma(t) = \left(x_0^i + \int_0^t A^i(\gamma(\tau))\,d\tau,\
t + \int_0^t B(\gamma(\tau))\,d\tau\right).
\end{displaymath}
It is
straightforward to check that
\begin{align*}
\int_0^t |A^i(\gamma(\tau))|\,d\tau
&\le \int_0^t C_0y(\gamma(t))^{\lambda-1}\,d\tau\\
&\le \int_0^t C_0C\tau^{\lambda-1}\,d\tau
\le C't^{\lambda}
\le C'\epsilon^{\lambda-\alpha}t^{\alpha};\\
\int_0^t |B(\gamma(\tau))|\,d\tau
&\le \int_0^t C_0y(\gamma(t))^{\lambda}\,d\tau\\
&\le \int_0^t C_0C\tau^{\lambda}\,d\tau
\le C't^{\lambda+1}
\le C'\epsilon^{\lambda-\alpha}t^{\alpha+1}.
\end{align*}
Thus if we
choose $\epsilon$
small enough, it follows that
$T$ maps $\scr X$ into $\scr X$.

We will show that, after choosing $\epsilon$ even smaller if necessary,
$T$ is a contraction.
Suppose $\gamma,\tw\gamma\in \scr X$.
If $(x^*,y^*)$ is any point along the line segment
between $\gamma(t)$ and $\tw\gamma(t)$,
\eqref{eq:b-gt-half-t}
implies that $\half t \le y^* \le \tfrac32 t$.
To estimate $d(T\gamma,T\tw\gamma)$, we use the mean-value theorem to obtain
\begin{align*}
t^{-\alpha}&\int_0^t \bigl|A^i(\gamma(\tau)) - A^i(\tw\gamma(\tau))\bigr|\,d\tau\\
&\le t^{-\alpha} \sum_j\int_0^t |\del_{x^j} A^i(x^*,y^*)|\, |a^j(\tau) - \tw a^j(\tau)|\,d\tau\\
&\qquad+ t^{-\alpha}\int_0^t |\del_y A^i(x^*,y^*)|\, |b(\tau) - \tw b(\tau)|\,d\tau \\
&\le t^{-\alpha}\int_0^t  mC_0(y^*)^{\lambda-1}\tau^\alpha d(\gamma,\tw\gamma)\,d\tau\\
&\qquad+ t^{-\alpha}\int_0^t C_0(y^*)^{\lambda-2} \tau^{\alpha+1}d(\gamma,\tw\gamma)\,d\tau \\
&\le t^{-\alpha} \int_0^t C_0C \tau^{\lambda+\alpha-1} d(\gamma,\tw\gamma)\,d\tau\\
&\le C' \epsilon^{\lambda}d(\gamma,\tw\gamma),
\end{align*}
for some $(x^*,y^*)$ on the line segment between $\gamma(\tau)$ and
$\tw\gamma(\tau)$.
The analogous estimate for the $y$-component is similar.
Thus $T$ is a contraction if we choose $\epsilon$ sufficiently small.  It therefore
has a unique fixed point in $\scr X$, which is a solution to
\eqref{eq:gamma-ode}.
By compactness, it is clear that for any given $C_1$ we can choose $\epsilon$
uniformly for $x_0\in K$.

To see that the solution is unique,
suppose $\gamma$ is any continuous solution to \eqref{eq:gamma-ode}. 
If we write $\gamma$ in the form \eqref{eq:gamma-a-b}, 
the equation $T\gamma=\gamma$ 
together with \eqref{eq:AB-estimates} implies successively that 
$|b(t)|\le Ct$, then 
$|b(t)| \le C t^{\lambda+1}$, and finally
$|a^i(t)| \le C t^\lambda$.
It follows that \eqref{eq:a-estimate} and 
\eqref{eq:b-estimate} hold
on $[0,\epsilon]$ if $C_1$ and $\epsilon$ are chosen appropriately.
Thus the restriction 
of $\gamma$ to $[0,\epsilon]$
is in $\scr X$, and so is equal 
to the unique fixed point of $T$.

Finally, we will prove the polyhomogeneity of the solution when
$V$ is polyhomogeneous.
Recall the spaces $\scr C^\delta$ defined 
at the beginning of Section \ref{section:polyhomogeneity}.
We will need the following fact about these spaces, which is proved by 
a straightforward analysis of the Taylor expansion of $F(u+f,v+g)$ about
$(u,v)=(u^1,\dots,u^m,v)$:
\begin{equation}\label{eq:Cdelta-composition}
\begin{gathered}
F\bigl(u(x,y)+f(x,y), v(x,y)+g(x,y)\bigr)
- F\bigl(u(x,y), v(x,y)\bigr)\\ 
\in \scr C^{\delta+\gamma}\\
\text{when }
F\in \scr C^\delta,\ 
u^i\in \scr C^0,\  
f^i\in \scr C^\gamma,\ 
v\in \scr C^1,\ 
g\in \scr C^{\gamma+1},\ 
\gamma\ge 0.
\end{gathered}
\end{equation}

For this proof, let $\scr A$ denote the space of polyhomogeneous
functions on $U$, and
for any $\delta\in \R$, define $\scr A^\delta = \scr A\cap \scr C^\delta$.
Thus a polyhomogeneous function $f$ is in 
$\scr A^\delta$ if and only if its leading term in the expansion
\eqref{eq:def-polyhomo} satisfies $s_1\ge\delta$ 
and, if $s_1=\delta$, $q_1=0$.
Our hypotheses imply that $A^i\in \scr A^{\lambda-1}$ and $B\in \scr A^{\lambda}$.

For each $x_0\in K$, 
let $a(x_0,t)$ and $b(x_0,t)$ denote the functions $a(t)$ and
$b(t)$ obtained above with initial condition $(x_0,0)$.
The standard argument showing that
solutions of ODEs depend smoothly upon initial values can be used
to obtain estimates of the following form for all multi-indices $\alpha$:
\begin{displaymath}
|\partial^\alpha_x a^i (x,t)|\le C_\alpha t^\lambda,
\quad |\partial^\alpha_x b (x,t)|\le C_\alpha t^{\lambda+1}.
\end{displaymath}
It is then straightforward to use  the differential equation to 
obtain estimates on
$(t\partial_t)^r\partial_x^\alpha a^i(x,t)$ and 
$(t\partial_t)^r\partial_x^\alpha b(x,t)$, showing that 
$a^i\in \scr C^\lambda$ and $b\in \scr C^{\lambda+1}$.

Suppose that for some integer $m\ge 1$, we have 
a ``partial polyhomogeneous expansion'' of the form
\begin{equation}\label{eq:partial-homogeneous}
\begin{gathered}
a^i(x,t) = p^i(x,t) + r^i(x,t), \qquad b(x,t) = q(x,t) + s(x,t),\\
\text{with }
p^i\in \scr A^\lambda,\ r^i\in \scr C^{m\lambda},\ 
q\in \scr A^{\lambda+1},\ s\in \scr C^{m\lambda+1}.
\end{gathered}
\end{equation}
The discussion in the preceding paragraph shows 
that \eqref{eq:partial-homogeneous} holds
with $m=1$ and  $p^i=q=0$.
Inserting \eq{eq:partial-homogeneous} into $T\gamma=\gamma$, 
we obtain $r^i = r^i_0 + r^i_1$, where 
\begin{align*}
r^i_0
&= \int_0^t A^i\bigl(x+p(x,\tau), \tau + q(x,\tau)\bigr)\,d\tau - p^i(x,t),\\
r^i_1 
&= \int_0^t \Bigl(A^i\bigl(x+p(x,\tau)+r(x,\tau), \tau+q(x,\tau)+s(x,\tau)\bigr)\\
&\qquad - A^i\bigl(x+p(x,\tau), \tau+q(x,\tau)\bigr)\Bigr)\,d\tau.
\end{align*}
From \eqref{eq:Cdelta-composition} we conclude that 
$r^i_1 \in\scr C^{(m+1)\lambda}$.  On the other hand, 
it is easy to check that $r^i_0$ is polyhomogeneous, and since
$r^i_0= r^i - r^i_1\in \scr C^{m\lambda}$, we have 
$r^i_0\in \scr A^{m\lambda}$.
A similar argument shows that $s = s_0+s_1$, with 
$s_0\in \scr A^{m\lambda+1}$ and $s_1\in\scr C^{(m+1)\lambda+1}$.
We let $P^i$ denote the sum of the 
(finitely many) terms in the expansion of 
$r^i_0$ that are in $\scr A^{m\lambda}\smallsetminus
\scr C^{(m+1)\lambda}$, and $Q$ the sum of the terms
in $s_0$ that are in $\scr A^{m\lambda+1}\smallsetminus
\scr C^{(m+1)\lambda+1}$.  Replacing $p^i$ by $p^i + P^i$ and $q$ by 
$q+ Q$,
we obtain \eqref{eq:partial-homogeneous} with $m+1$ in place of 
$m$.  Continuing by induction, we conclude that $a$ and $b$ are polyhomogeneous.
Since both are in $\scr A^\lambda$, it follows that 
$(x,t)\mapsto\gamma_x(t)$ is a 
$C^{0,\lambda}$ map.
\end{proof}

We are finally ready to prove the main theorem.

\begin{proof}[Proof of Theorem \ref{thm:main}]\
Suppose $g$ satisfies the hypotheses of the theorem.
Let $\Psi$ be the collar diffeomorphism given by
Theorem \ref{thm:harmonic-map-existence}, and
let $\tw g = \Psi^* g$, which is polyhomogeneous
and conformally compact of class $C^{1,\lambda}$, $0<\lambda<1$.
Theorem 
\ref{thm:harmonic-map-existence} shows that 
$\rho^2\Psi^* g - \rho^2 h\in C^{1,0}_{3+\lambda}(M_R;\Sigma^2)$,
which implies that $\rho^2\Psi^* g - \rho^2 h = O(\rho^{1+\lambda})$ in 
background coordinates, so $\Psi^* g$ and $h$ have the
same conformal infinity $\tw\gamma$.

Then let $\Gamma$ be the collar diffeomorphism
given by Lemma \ref{lemma:phg-Fermi-coords}, so that $\Gamma^*\tw g$
has the form
\eqref{eq:metric-in-Fermi-coords}.
Because 
$r/\rho\to 1$ at $\del M$, it follows that $\Gamma^*\tw g$ 
also has $\tw \gamma$ as conformal infinity.
Because 
$\rho^2
\Gamma^*\tw g$ is continuous up to $\del M$ and has a smooth
restriction to $\del M$, it follows that the log terms in the
asymptotic expansion for 
$G(\rho)$ all occur with positive powers
of $\rho$. Once polyhomogeneity is known, the detailed form of the
expansion is established by matching powers and coefficients of
various terms appearing in the equations. Such a study has been
done by Robin Graham and Charles Fefferman \cite{FG2003} (see also
\cite{KichenassamyFG}), and the results there imply that the
expansions are of the form described in Theorem
\ref{thm:main}.
In the special case when 
$\dim M=3$, it is possible to choose the 
conformal infinity $\tw\gamma$ to have constant Gaussian curvature, and
an easy computation shows that the first log term
vanishes in this case
(cf. ~ \cite{KichenassamyFG,FG,FG2003}),
so $\rho^2 \Gamma^*\tw g$ is always smooth.

The proof is completed by letting $\Phi$ be the collar diffeomorphism
$\Psi\circ \Gamma$.
\end{proof}

\bibliographystyle{amsplain}
\bibliography{bregularity}

\def\cprime{$'$}
\providecommand{\bysame}{\leavevmode\hbox to3em{\hrulefill}\thinspace}
\providecommand{\MR}{\relax\ifhmode\unskip\space\fi MR }
\providecommand{\MRhref}[2]{%
  \href{http://www.ams.org/mathscinet-getitem?mr=#1}{#2}
}
\providecommand{\href}[2]{#2}
\begin{thebibliography}{10}

\bibitem{mand2}
M.~T. Anderson, \emph{{Einstein} metrics with prescribed conformal infinity on
  $4$-manifolds}, math.DG/0105243, 2001/2004.

\bibitem{mand1}
\bysame, \emph{Boundary regularity, uniqueness and non-uniqueness for {AH
  Einstein} metrics on $4$-manifolds}, Adv.\ Math. \textbf{179} (2003),
  205--249, MR2010802, Zbl 1048.53032.

\bibitem{AndChDiss}
L.~Andersson and P.~T. Chru\'{s}ciel, \emph{Solutions of the constraint
  equations in general relativity satisfying ``hyperboloidal boundary
  conditions''}, Dissertationes Math. \textbf{355} (1996), 1--100, MR1405962,
  Zbl 0873.35101.

\bibitem{Besse}
A.~L. Besse, \emph{{Einstein manifolds}}, Ergeb.\ Math.\ Grenzgeb.\ (3),
  vol.~10, Springer Verlag, Berlin, New York, Heidelberg, 1987, MR0867684, Zbl
  0613.53001.

\bibitem{biquard}
O.~Biquard, \emph{{M\'etriques d'Einstein asymptotiquement sym\'etriques
  (Asymptotically symmetric Einstein metrics)}}, {Ast\'erisque 265, Paris:
  Soci\'et\'e Math\'ematique de France, 109 pp.}, 2000, MR1760319, Zbl
  0967.53030.

\bibitem{ChMS}
P.~T. Chru{\'s}ciel, M.~A.~H. MacCallum, and D.~Singleton, \emph{Gravitational
  waves in general relativity. {XIV}: {B}ondi expansions and the
  ``polyhomogeneity'' of {S}cri}, Philos.\ Trans.\ Roy.\ Soc.\ London Ser.\ A
  \textbf{350} (1995), 113--141, MR1325206, Zbl 0829.53065.

\bibitem{FG}
C.~Fefferman and C.~R. Graham, \emph{{Conformal invariants}}, {\'Elie Cartan et
  les math\'ematiques d'aujourd'hui, The mathematical heritage of \'Elie
  Cartan, S\'emin. Lyon 1984, Ast\'erisque, Num\'ero Hors S\'er., 95--116},
  1985, MR0837196, Zbl 0602.53007.

\bibitem{FG2003}
\bysame, \emph{{\it The ambient metric}}, in preparation, 2005.

\bibitem{GL}
C.~R. Graham and J.~M. Lee, \emph{Einstein metrics with prescribed conformal
  infinity on the ball}, Adv.\ Math. \textbf{87} (1991), 186--225, MR1112625,
  Zbl 0765.53034.

\bibitem{deHaro:2000xn}
{S.~de} Haro, S.~N. Solodukhin, and K.~Skenderis, \emph{Holographic
  reconstruction of spacetime and renormalization in the {AdS/CFT}
  correspondence}, Comm.\ Math.\ Phys. \textbf{217} (2001), 595--622,
  hep-th/0002230, MR1822109, Zbl 0984.83043.

\bibitem{HG}
K.~Hirachi and C.~R. Graham, \emph{The ambient obstruction tensor and
  {$Q$}-curvature}, math.DG/0405068, 2004.

\bibitem{KichenassamyFG}
Satyanad Kichenassamy, \emph{On a conjecture of {F}efferman and {G}raham},
  Adv.\ Math. \textbf{184} (2004), no.~2, 268--288, MR2054017, Zbl pre02081968.

\bibitem{Lee:spectrum}
J.~M. Lee, \emph{The spectrum of an asymptotically hyperbolic {E}instein
  manifold}, Comm.\ Anal.\ Geom. \textbf{3} (1995), 253--271, MR1362652, Zbl
  0934.58029.

\bibitem{Lee:fredholm}
\bysame, \emph{Fredholm operators and {E}instein metrics on conformally compact
  manifolds}, math.DG/0105046, to appear in Mem.\ Amer.\ Math.\ Soc., 2005.

\bibitem{llave-obaya}
{R.~de la} Llave and R.~Obaya, \emph{Regularity of the composition operator in
  spaces of {H}\"older functions}, Discrete Contin. Dyn. Syst. \textbf{5}
  (1999), 157--184, MR1664481, Zbl 0956.47029.

\bibitem{mazzeo:hodge}
R.~Mazzeo, \emph{{The Hodge cohomology of a conformally compact metric.}}, J.\
  Differential Geom. \textbf{28} (1988), 309--339, MR0961517, Zbl 0656.53042.

\bibitem{mazzeo}
\bysame, \emph{Elliptic theory of differential edge operators {I}}, Comm.
  Partial Differential Equations \textbf{16} (1991), 1615--1664, MR1133743, Zbl
  0745.58045.

\end{thebibliography}

\end{document}